\documentclass{amsart}

\usepackage{amsfonts,amsmath,amssymb}
\usepackage{graphicx}
\newtheorem{theorem}{Theorem}
\newtheorem{lemma}[theorem]{Lemma}
\newtheorem{definition}[theorem]{Definition}
\newtheorem{example}[theorem]{Example}

\newtheorem{remark}[theorem]{Remark}

\title{Desingularization in Computational Applications and Experiments}
\author{Anne Fr\"uhbis-Kr\"uger}
\address{Institut f\"ur Algebraische Geometrie, Leibniz Universit\"at Hannover,
         Germany}
\begin{document}
\thanks{partially supported by DFG priority program SPP 1489 'Algorithmic and
        Experimental Methods in Algebra, Geometry and Number Theory'}

\maketitle

\begin{abstract}
After briefly recalling some computational aspects of blowing up
and of representation of resolution data common to a wide range of
desingularization algorithms (in the general case as well as
in special cases like surfaces or binomial varieties), we shall
proceed to computational applications of resolution of singularities
in singularity theory and algebraic geometry, also touching on 
relations to algebraic statistics and machine learning.
Namely, we explain how to compute the intersection form and dual graph of 
resolution for surfaces, how to determine discrepancies, the log-canoncial 
threshold and the topological Zeta-function on the basis of desingularization
data. We shall also briefly see
how resolution data comes into play for Bernstein-Sato polynomials, and we
mention some settings in which desingularization algorithms can be used 
for computational experiments. The latter is simply an invitation
to the readers to think themselves about experiments using existing software,
whenever it seems suitable for their own work. 
\end{abstract}

\section{Introduction}

This article originated from the notes of an invited talk at the Clay 
Mathematics Institute summer school on "The Resolution of Singular Algebraic 
Varieties" in Obergurgl, Austria, 2012. As the whole school was devoted to 
desingularization, the focus in this particular contribution is on 
applications and on practical aspects. A general knowledge of resolution of
singularities and different approaches to this task is assumed and can be
acquired from other parts of this proceedings volume. The overall goal
of this article is to give readers a first impression of a small choice of 
applications and point them to good sources for further reading on each
of the subjects. A detailed treatment of each of the topics would fill an 
article by itself and is thus beyond the scope here.\\
One focus here is on the practical side. To this end, we first 
revisit desingularization algorithms in section 2 and have a closer look
at the representation of resolution data: as a consequence of the heavy use
of blowing up, the data is distributed over a rather large number of charts
which need to be glued appropriately. Glueing, however, only describes the
theoretical side of the process; from the practical point of view, it is 
closer to an identification of common points in charts. \\
In section 3 we focus on applications needing different amounts of 
resolution data. Using an abstract resolution of singularities only, the
computation of the intersection form and dual graph of the resolution for
surface singularities requires the smallest amount of data. For determining
discrepancies and the log-canonical threshold, we already need an embedded
resolution which is also required for the third application, the computation
of the topological zeta function.\\
In the rather short last section, we only sketch two settings in which one
might want to use algorithmic desingularization as an experimental tool:
the roots of the Bernstein-Sato polynomial and resolution experiments in
positive characteristic. This last part is not intended to provide actual
research projects. It is only intended to help develop a feeling for settings 
in which experiments can be helpful.\\
All parts of the article are illustrated by the same simple 
example which is desingularized using a variant of Villamayor's algorithm
available in {\sc Singular}. Based on this resolution data, all further
applications are also accompanied by the corresponding {\sc Singular}-code.
The {\sc Singular} code is not explained in detail, but hints and explanations
on the appearing commands and their output are provided as comments in the 
examples. To further familiarize with the use of {\sc Singular} in this 
context, we recommend that the readers try out the given session themselves 
and use the built-in manual of {\sc Singular} to obtain further information 
on the commands (e.g.: {\tt help resolve;} returns the help page of the 
command {\tt resolve}).\\
I would like to thank the organizers of the summer school for the invitation.
Insights from conversations with many colleagues have contributed to the 
content of these notes. In particular, I would like to thank Herwig Hauser,
Gerhard Pfister, Ignacio Luengo, Alejandro Melle, Frank-Olaf Schreyer,
Wolfram Decker, Hans Sch\"onemann, Duco van Straten, Nobuki Takayama, 
Shaowei Lin, Frank Kiraly, Bernd Sturmfels, Zach Teitler, Nero Budur, 
Rocio Blanco and Santiago Encinas for comments, of which some led to 
applications explained here and some others helped seeing the applications 
in broader context. 
For reading earlier versions of this article and many helpful questions on 
the subject of this article, I 
am indepted to Frithjof Schulze and Bas Heijne.

\section{Desingularization from the computational side}

Before turning toward applications of resolution of singularities, we need
to review certain aspects of algorithmic desingularization to understand
the way in which the computed resolution data is represented. Although there
are various settings in which different resolution algorithms have been
created, we may discern three main approaches suitable for the purposes of
this article: the algorithms based on Hironaka's proof in characteristic 
zero in any dimension (see e.g. \cite{BEV}, \cite{BM}, \cite{EH}), the 
algorithms for binomial and toric ideals (see e.g. \cite{Ful}, \cite{Bl}, 
\cite{BM2}) and the algorithms for 2-dimensional varieties and schemes (such 
as \cite{Beck} -- based on \cite{Jung} -- or \cite{Lip}).\footnote{Of 
course this list of approaches is far from exhaustive, but it is intended 
to narrow down our scope to those which lead to similar forms of resolution 
data allowing similar applications later on.}
The algorithms of the first kind of approach involve embedded 
desingularization, i.e. they blow up a smooth ambient space and consider the 
strict transform of the variety and the exceptional divisors inside the new 
ambient space. The algorithms for 2-dimensional varieties on the other hand, 
do not consider the embedded situation, but blow up the variety itself and 
consider exceptional divisors inside the blown up variety.\\

It would be beyond the scope of this article to cover all these algorithms
in depth, but they all have certain ingredients in common.
In the first two situations, a desingularization is achieved by finite 
sequences of blow-ups at suitable non-singular centers. The differences 
between these algorithms then lie in the choice of center, which is the 
key step of each of these, but does not affect the structure of the practical 
representation of resolution data. For the third class of algorithms, 
blow-ups are not the only tool, but are combined with other tools, in 
particular normalization steps. However, the exceptional divisors to be 
studied arise from blow-ups and additionally only require proper tracing 
through the normalization steps if necessary. Therefore the technique to 
focus on in this context is blowing up; more precisely blowing up at 
non-singular centers.\\

\subsection{Blowing up -- the computational side}

Let us briefly recall the definition of blowing up, as it can be found in 
any textbook on algebraic geometry (e.g. \cite{Har}), before explaining its 
computational side:

\begin{definition}
Let $X$ be a scheme and $Z \subset X$ a subscheme corresponding to a coherent
ideal sheaf ${\mathcal I}$. The blowing up of $X$ with center $Z$ is 
$$\pi:\overline{X}:=Proj(\bigoplus_{d\geq0} {\mathcal I}^d) \longrightarrow X.$$ 
Let $Y {\stackrel{i}{\hookrightarrow}} X$ be a closed subscheme and 
$\pi_1:\overline{Y} \longrightarrow Y$ the blow up of $Y$ along 
$i^{-1}{\mathcal I}{\mathcal O}_Y$. Then the following diagram commutes 
\begin{eqnarray*}
\overline{Y} & \hookrightarrow & \overline{X}\\
\pi_1 \downarrow & & \downarrow \pi\\
Y & \hookrightarrow & X
\end{eqnarray*}
$\overline{Y}$ is called the strict transform of $Y$, $\pi^*(Y)$ the total 
transform of $Y$.
\end{definition}

To make $\overline{X}$ accessible to explicit computations of examples in 
computer algebra systems, it should best be described as the zero set of an 
ideal in a suitable ring. We shall assume now for simplicity of
presentation that $X$ is affine because schemes are usually represented 
in computer algebra systems by means of affine covers. So we are dealing 
with the following situation: $J=\langle f_1,\dots,f_m\rangle \subset A$
is the vanishing ideal of the center $Z \subset X = Spec(A)$ and the task is 
to compute
$$Proj(\bigoplus_{d\geq 0} J^d).$$
To this end, we consider the canonical graded $A$-algebra homomorphism
$$\Phi: A[y_1,\dots,y_m] \longrightarrow 
  \bigoplus_{n\geq 0} J^nt^tn \subset A[t]$$
defined by $\Phi(y_i)=tf_i$. The desired object $\bigoplus_{d \geq 0} J^d$
is then isomorphic to 
$$A[y_1,\dots,y_m]/ker(\Phi)$$ 
or from the more geometric point of view 
$V(ker(\Phi)) \subset Spec(A) \times {\mathbb P}^{m-1}$.\\

The computation of the kernel in the above considerations is a standard
basis computation and further such computations arise during the calculation
of suitable centers in the different algorithms. Moreover, each blowing up 
introduces a further set of new variables as we have just seen and 
desingularization is hardly ever achieved with just one or two blow-ups --
usually we are seeing long sequences thereof.  The performance of standard
bases based algorithms, on the other hand, is very sensitive to the number 
of variables as its complexity is doubly exponential in this number. Therefore
it is vital from the practical point of view to pass from the 
${\mathbb P}^{m-1}$ to $m$ affine charts and pursue the resolution process
further in each of the charts\footnote{In practice, it is very useful to
discard all charts not appearing in any other chart and not containing any 
information which is relevant for the desired application.}; 
consistency of the choice of centers does not pose a problem at this stage 
as this follows from the underlying desingularization algorithm. Although 
this creates an often very large tree of charts with the final charts 
being the leaves and although it postpones a certain part of the work, this 
approach has a further important advantage: it allows parallel computation 
by treating several charts on different processors or cores at the same time.\\

As a sideremark, we also want to mention the computation of the transforms,
because they appeared in the definition cited above; for simplicity of 
notation we only state the affine case. Without any computational
effort, we obtain the exceptional divisor as 
$I(H)=J{\mathcal O}_{\overline{X}}$ and the total transform 
$\pi^*(I)=I{\mathcal O}_{\overline{X}}$ for a subvariety $V(I)$ in our affine
chart. The strict transform is then obtained by a saturation, i.e. an iteration
of ideal quotients until it stabilizes: 
$I_{\overline{V(I)}}=(\pi^*(I): I(H)^\infty)$; for the weak transform, the
iteration stops before it stabilizes, namely at the point where multiplying
with the ideal of the exceptional divisor gives back the result of the 
previous iteration step for the last time. These ideal quotient computations
are again based on standard bases.\\

\subsection{Identification of Points in Different Charts}

As we have just seen, it is more useful for the overall performance to pass 
to affine charts after each blow-up, even though this leads to an often 
rather large tree of charts. As a consequence, it is not possible to directly 
work with the result without any preparation steps, namely identifying points
which are present in more than one chart -- a practical step equivalent to the 
glueing of the charts. This is of particular interest for the identification of 
exceptional divisors which are present in more than one chart. 

For the identification of points in different charts, we need to pass through
the tree of charts -- from one final chart all the way back to the last
common ancestor of the two charts and then forward to the other final chart.
As blowing up is an isomorphism away from the center, this step does not 
pose any problems as long as we do not need to identify points lying on
an exceptional divisor which was not yet created in the last common ancestor
chart. In the latter case, however, we do not have a direct means of keeping
track of points originating from the same point of the center. The way out
of this dilemma is a representation of points on the exceptional divisor as
the intersection of the exceptional divisor with an auxilliary variety not
contained in the exceptional divisor. More formally, the following simple 
observation from commutative algebra can be used:

\begin{remark}
Let $I \subset K[x_1,\dots,x_n]$ be a prime ideal, $J \subset K[a_1,\dots,x_n]$
another ideal such that $I+J$ is equidimensional and $ht(I)= ht(I+J)- r$ for
some integer $0 < r< n$. Then there exist polynomials $p_1,\dots,p_r \in I+J$
and a polynomial $f \in K[x_1,\dots,x_n]$ such that
$$\sqrt{I+J} = \sqrt{(I+\langle p_1,\dots,p_r\rangle):f}.$$
\end{remark}

At first glance, this seems non-constructive, but it turns out to be
applicable in a very convenient way: In our situation, $I$ describes the 
intersection of exceptional divisors containing the point (or subvariety 
$V(J)$) in question. As any sufficiently general set of polynomials 
$p_1,\dots,p_r \in J \setminus (I \cap J)$ with correct height of 
$I+\langle p_1,\dots,p_r\rangle$ will serve our purpose and so will any
$f$ that excludes all extra components of $I + \langle p_1,\dots,p_r\rangle$,
we have sufficient freedom of choice of $p_1,\dots,p_r,f$ such that none of 
these is contained in any further exceptional divisor whose moment of birth 
has to be crossed during the blowing down process. \\

Given this means of identification of points, we can now also identify
the exceptional divisors or, more precisely, the centers leading to the
respective exceptional divisors. To avoid unnecessary comparisons between
centers in different charts, we can a priori rule out all comparisons 
involving centers lying in different exceptional divisors. If the 
desingularization is controlled by an invariant as in \cite{BEV}, \cite{BM}
or \cite{Bl}, we can also avoid comparisons with different values of 
the controlling invariant, because these cannot give birth to the same 
exceptional divisor either. 

\begin{example}
To illustrate the explanations given so far and to provide a practical example
to be used for all further applications, we now consider an isolated 
surface singularity of type $A_4$ at the origin. We shall illustrate this
example using the computer algebra system {\sc Singular} (\cite{Sing}).
\begin{verbatim}
> // load the appropriate libraries for resolution of singularities and
> // applications thereof
> LIB"resolve.lib";
> LIB"reszeta.lib";
> LIB"resgraph.lib";

> // define the singularity
> ring R=0,(x,y,z),dp;
> ideal I=x5+y2+z2;               // an A4 surface singularity

> // compute a resolution of the singularity (Villamayor-approach)
> list L=resolve(I);
> size(L[1]);                     // final charts
6
> size(L[2]);                     // all charts
11
> def r9=L[2][9];                 // go to chart 9
> setring r9;
> showBO(BO);                     // show data in chart 9
                       
==== Ambient Space: 
_[1]=0
      
==== Ideal of Variety: 
_[1]=x(2)^2+y(0)+1
      
==== Exceptional Divisors: 
[1]:
   _[1]=1
[2]:
   _[1]=y(0)
[3]:
   _[1]=1
[4]:
   _[1]=x(1)
   
==== Images of variables of original ring:
_[1]=x(1)^2*y(0)
_[2]=x(1)^5*x(2)*y(0)^2
_[3]=x(1)^5*y(0)^2
> setring R;                     // go back to old ring
\end{verbatim} 
This yields a total number of 11 charts of which 6 are final charts. All
blow-ups have zero-dimensional centers except the two line blow-ups leading
from chart 6 to charts 8 and 9 and from chart 7 to charts 10 and 11. As it
would not be very useful to reproduce all data of this resolution here, we show
the total tree of charts as illustration \ref{TreeA4} and give the content 
of chart 9 as an example:\\[0.3cm]
\begin{figure}[h]\includegraphics[width=6cm]{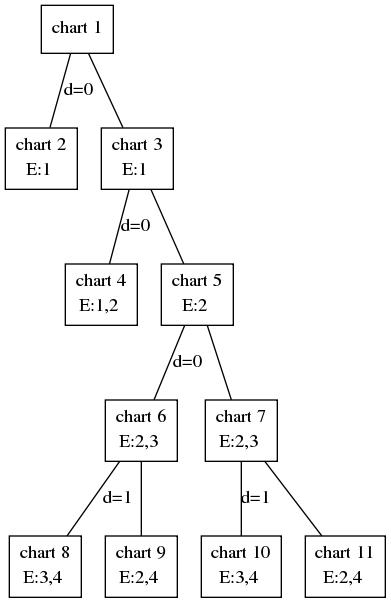} \label{TreeA4}
\caption{Tree of charts of an embedded desingularization of an $A_4$ surface
singularity. The numbers given in the second line in each chart are the labels
of the exceptional divisors visible in the respective chart. The numbers stated
as $d=0$ or $d=1$ state the dimension of the center of the corresponding 
blow-up. Charts providing only data which is also present in other charts are
not shown.}\end{figure}
\phantom{abc}strict transform\footnote{Weak and strict transform coincide for hypersurfaces.}:  $V(x_2^2+y_0+1)$\\
\phantom{abc}exceptional divisors: $V(y_0)$ from 2nd blow-up\\
\phantom{abcexceptional divisors:} $V(x_1)$ from last blow-up\\
\phantom{abc}images of variables of original ring:\\
\phantom{abcdefg} $x \longmapsto x_1^2y_0$\\
\phantom{abcdefg} $y \longmapsto x_1^5x_2y_0^2$\\
\phantom{abcdefg} $z \longmapsto x_1^5y_0^2$\\[0.3cm]
For example, the exceptional divisor originating from the last blow-up leading
to chart 9 needs to be compared to the exceptional divisor originating from
the last blow-up leading to chart 10. To this end, one needs to consider the 
centers computed in charts 6 and 7 which turn out to be the intersection of 
the exceptional divisors labeled $2$ and $3$, if one considers the output
of the resolution process in detail\footnote{We encourage the reader to 
verify this by typing the above sequence of commands into {\sc Singular} 
and then exploring the data in the different charts.}. 
Therefore the last exceptional divisors in charts 9 and 10 coincide. This 
identification is implemented in {\sc Singular} and can be used in the 
following way: 
\begin{verbatim}
> // identify the exceptional divisors
> list coll=collectDiv(L); 
> coll[1];
0,0,0,0,  // no exc. div. in chart 1
1,0,0,0,  // first exc. div. is first in chart 2
1,0,0,0,  // ..... 
1,2,0,0,  // data too hard to read, better
0,2,0,0,  // use command below to create figure 1
0,2,3,0,
0,2,3,0,
0,0,3,4,
0,2,0,4,
0,0,3,4,
0,2,0,4 

> //present the tree of charts as shown in figure 1
> ResTree(L,coll[1]);
\end{verbatim}
\end{example}

\section{Applications of Resolution of Singularities}

The applications we present in this section originate from different 
subfields of mathematics ranging from algebraic geometry to singularity 
theory and $D$-modules. For each application we shall revisit our example 
from the previous section and also show how to perform the corresponding
computation using {\sc Singular}.
  
\subsection{Intersection Form and Dual Graph of Resolution}

Given a resolution of an isolated surface singularity, we want to compute
the intersection matrix of the excpetional divisors. This task does not 
require an embedded resolution of singularities, only an abstract one. Given
such a desingularization, it can then be split up into 3 different subtasks:

\begin{itemize}
\item[(1)] computation of the intersections $E_i . E_j$ for exceptional
           curves $E_i \neq E_j$
\item[(2)] computation of the self-intersection numbers $E_i^2$ for the
           exceptional curves $E_i$
\item[(3)] representation of the result as the dual graph of the resolution
\end{itemize}

If the given resolution, is not an abstract one, but an embedded one - like
the result of Villamayor's algorithm - we need to add a preliminary step
\begin{itemize}
\item[(0)] determine an abstract resolution from an embedded one.
           \footnote{This 
           is achieved by canceling all trailing blow-ups in our tree of charts
           which are unnecessary for the non-embedded case. Then the 
           intersection of the remaining excpetional divisors with the strict
           transform yields the exceptional locus of the non-embedded 
           resolution.}
\end{itemize}

Although the definition of intersection numbers of divisors on surfaces
can be found in many textbooks on algebraic geometry (e.g. in \cite{Har}, V.1),
we give a brief summary of the used properties for readers' convenience:

\begin{definition}
Let $D_1,D_2$ be divisors in general position\footnote{$D_1$ and $D_2$ are
in general position, if the intersection $Supp(D_1) \cap Supp(D_2)$ is either
empty or a finite set of points.} on a non-singular surface $X$. Then the
intersection number of $D_1$ and $D_2$ is defined as
$$D_1 . D_2 := \sum_{x \in D_1 \cap D_2} (D_1 . D_2)_x$$
where $(D_1 . D_2)_x$ denotes the intersection multiplicity of $D_1$ and 
$D_2$ at $x$.
\end{definition}

\begin{lemma}
For any divisors $D_1$ and $D_2$ on a non-singular surface $X$, there exist
divisors $D_1'$ and $D_2'$, linearly equivalent to $D_1$ and $D_2$ 
respectively, such that $D_1'$ and $D_2'$ are in general position. 
\end{lemma}

\begin{lemma}
Intersection numbers have the following basic properties:
\begin{enumerate}
\item[(a)] For any divisors $C$ and $D$: $C.D = D.C$.
\item[(b)] For any divisors $C$, $D_1$ and $D_2$: $C.(D_1+D_2) = C.D_1 + C.D_2$
\item[(c)] For any divisors $C$, $D_1$ and $D_2$, such that $D_1$ and
           $D_2$ are linearly equivalent: $C.D_1 = C.D_2$.
\end{enumerate}
\end{lemma}

At this point, we know what we want to compute, but we have to take care of
another practical problem before proceeding to the actual computation which
is then a straight forward calculation of the intersection numbers of  
exceptional curves $E_i \neq E_j$. 
The practical problem is that computations in a 
computer algebra system usually take place in polynomial rings over the
rationals or algebraic extensions thereof. So we easily achieve a decomposition
of the exceptional divisor into ${\mathbb Q}$-irreducible components, but 
we need to consider ${\mathbb C}$-irreducible components to obtain the
intersection matrix we expect from the theoretical point of view. To this end,
passing to suitable extensions of the ground field may be necessary -- be it
explicitly by introducing a new variable and a minimal polynomial (slowing 
down subsequent computations) or implicitly by taking into account the number
of components over ${\mathbb C}$ for each ${\mathbb Q}$-component.  

\begin{example}
Revisiting our example of a desingularization of an $A_4$ surface singularity,
we first need to pass to an abstract resolution. 
\begin{verbatim}
\\ compute part of the tree of charts relevant for abstract resolution
> abstractR(L)[1];
   0,1,0,1,1,0,0,0,0,0,0   //final charts are 2,4,5
> abstractR(L)[2];
   0,0,0,0,0,1,1,1,1,1,1   //charts 6 and higher are irrelevant 
                           //for non-embedded case
\end{verbatim}
So we only see 2 exceptional divisors in the final charts, the ones labeled
1 and 2 (cf. figure \ref{TreeA4}). But looking at the charts in more detail, 
we would see e.g. in chart 4 that the first divisor is given by 
$V(y_2,y_0^2+1)$ and the second one by $V(x_2,y_0^2+1)$. So each of these 
has two ${\mathbb C}$-irreducible components. Considering these 
${\mathbb C}$-components, we can obtain the following intersection data 
(seeing two of the intersections directly in chart 4 and the remaining one 
in chart 5):
$$\left(\begin{array}{cccc}
           * & 0 & 1 & 0\cr
           0 & * & 0 & 1\cr
           1 & 0 & * & 1\cr
           0 & 1 & 1 & *\end{array}\right)$$
Hence only the self-intersection numbers -- marked as * in the matrix --
are still missing.
\end{example}

For the self-intersection numbers of the exceptional curves, we need to
make use of another property of divisors in the context of desingularization:

\begin{lemma}
Let $\pi: \tilde{X} \longrightarrow X$ be a resolution of singularities of a
surface $X$. Let $D_1$ be a divisor on $\tilde{X}$ all of whose components
are exceptional curves of $\pi$ and let $D_2$ be a divisor on $X$, then
$$\pi^*(D_2).D_1 =0.$$
\end{lemma}

Denoting by $E_1,\dots,E_s$ the ${\mathbb C}$-irreducible exceptional curves,
we can hence consider a linear form $h: X \longrightarrow {\mathbb C}$
passing through the only singular point of $X$ and the divisor $D$ defined
by it. We then know 
$$\pi^*(D)=\sum_{i=1}^s c_iE_i + H$$
where $H$ denotes the strict transform of $D$ and the $c_i$ are suitable 
integers. From the lemma we additionally know that
$$0 = \pi^*(D).E_j = \sum_{i=1}^s c_iE_i.E_j + H.E_j\;\; 
      \forall 1 \leq j \leq s$$
where all intersection numbers are known or directly computable in each
of the equations except the self-intersection numbers $E_j.E_j$ which we
can compute in this way. For the dual graph of the resolution, each 
divisor is represented by a vertex (those with self-intersection -2 are
unlabeled, the other ones labeled by their self-intersection number),
each intersection is represented by an edge linking the two vertices 
corresponding to the intersecting exceptional curves.

\begin{example}
The computation of the intersection form is implemented in {\sc Singular}
and can be used as follows:
\begin{verbatim}
> // intersection matrix of exceptional curves
> // (no previous abstractR is needed, this is done automatically 
> // in the procedure intersectionDiv)
> list iD=intersectionDiv(L);
> iD[1];  
-2,0,1,0,
0,-2,0,1,
1,0,-2,1,
0,1,1,-2 

> // draw the dual graph of the resolution
InterDiv(iD[1]);  
\end{verbatim}
This yields the expected intersection matrix (as entry iD[1] of the result)
$$\left(\begin{array}{rrrr}
          -2 & 0 & 1 & 0\cr
           0 &-2 & 0 & 1\cr
           1 & 0 &-2 & 1\cr
           0 & 1 & 1 &-2\end{array}\right)$$
and the dual graph of the resolution which is just the Dynkin diagram of 
the $A_4$ singularity.
\end{example}

\subsection{Discrepancies and Log-canonical Threshold}

In contrast to the last task, which only required an abstract resolution of 
the given surface singularity, the task of computing (log-)discrepancies
and the log-canonical threshold requires embedded desingularization (or
principalization of ideals). This is provided by Villamayor's algorithm.
As before we first recall the definitions and some basic properties (see e.g.
\cite{Must} 
for a direct and accessible introduction to the topic). To keep the exposition
as short as possible and the considerations directly accesible to explicit 
computation, we restrict our treatment here to the case of a singular affine
variety.  

\begin{definition}
Let $f \in {\mathbb C}[x_1,\ldots,x_n]$ be a non-zero polynomial defining a
hypersurface $V$ and let $\pi :X \longrightarrow {\mathbb C}^n$ be an embedded
resolution of $V$. Denote by $E_i$, $i \in I$, the irreducible components
of the divisor $\pi^{-1}(f^{-1}(0))$. Let $N(E_j)$ denote the multiplicity 
of $E_j$, $j \in J$ in the divisor of $f \circ \pi$ and let $\nu(E_j)-1$ be
the multiplicity of $E_j$ in the divisor 
$K_{X/{\mathbb C}^n}=\pi^*(dx_1\wedge\ldots\wedge dx_n)$. Then the 
log-discrepancies of the pair $({\mathbb C}^n,V)$ w.r.t. $E_j$, $j \in J$,
are
$$a(E_j;{\mathbb C}^n,V) := \nu(E_j) - N(E_j).$$
The minimal log discrepancy of the pair $({\mathbb C}^n,V)$ along a 
closed subset $W \subset {\mathbb C}^n$ is the minimum over the 
log-discrepancies for all $E_j$ with $\pi(E_j) \subset W$, i.e. originating 
from (sequences of) blow-ups with centers in $W$.
The log-canonical threshold of the pair $({\mathbb C}^n,V)$ is defined as
$$ lct({\mathbb C}^n,V) = inf_{j \in J} \frac{\nu(E_j)}{N(E_j)}.$$
\end{definition}

\begin{remark}
The above definition of log-discrepancies and log-canonical threshold holds
in a far broader context. Allowing more general pairs $(Y,V)$ it is also the
basis for calling a resolution of singularities log-canonical, if the minimal
log discrepancy of the pair along all of $Y$ is non-negative, and log-terminal,
if it is positive.
\end{remark}

As we already achieved an identification of exceptional divisors in a previous 
section, the only computational task here is the computation of the 
multiplicities $N(E_i)$ and $\nu(E_i)$. The fact that we might be dealing 
with ${\mathbb Q}$-irreducible, but ${\mathbb C}$-reducible $E_i$ does not
pose any problem here, because we can easily check that the respective 
multiplicities coincide for all ${\mathbb C}$-components of the same 
${\mathbb Q}$-component. To compute these multiplicities from the resolution
data in the final charts (i.e. without moving through the tree of charts),
we can determine $N(E_i)$ by finding the highest exponent $j$ such that
$I(E_j)^j:J$ is still the whole ring where $J$ denotes the ideal of the
total transform of the original variety. To compute the $\nu(E_i)$ we can
use the same approach, but taking into account the appropriate Jacobian 
determinant.

\begin{example}
We now simply continue our {\sc Singular} session on the basis of the data 
already computed in the previous examples
\begin{verbatim}
>// identify exceptional divisors (embedded case)
> list iden=prepEmbDiv(L);

>// multiplicities N(Ei)
> intvec cN=computeN(L,iden);
> cN;           // last integer is strict transform
2,4,5,10,1

>// multiplicities v(Ei)
> intvec cV=computeV(L,iden);
> cV;           // last integer is strict transform
3,5,7,12,1

>// log-discrepancies
> discrepancy(L);
0,0,1,1

>// compute log-canonical threshold
>// as an example of a loop in Singular
> number lct=number(cV[1])/number(cN[1]);
> number lcttemp;
> for (int i=1; i < size(cV); i++)
> {
>   lcttemp=number(cV[i])/number(cN[i]);
>   if(lcttemp < lct) 
>   {
>      lct=lcttemp;
>   }
> }
> lct;
6/5
\end{verbatim}
\end{example}

The log-canonical threshold, in particular, is a very important invariant
which appears in many different contexts, ranging from a rather direct 
study of properties of pairs to the study of multiplier ideals (cf. 
\cite{Laz}), to motivic integration or to Tian's $\alpha$-invariant which 
provides a criterion for the existence of K\"ahler-Einstein metrics (cf. 
\cite{Tian}). \\

A real analogue to the log-canonical threshold, the real log-canonical 
threshold appears when applying resolution of singularities in the real 
setting \cite{Sai}. In algebraic statistics, more precisely in model 
selection in Bayesian statistics, desingularization plays an important 
role in understanding singular models by monomializing the so-called 
Kullback-Leibler function at the true distribution. In this context the 
real log-canonical threshold is then used to study the asymptotics of the 
likelihood integral.\cite{Lin}. It also appears in singular learning 
theory as the learning coefficient.

\subsection{Topological Zeta Function}

Building upon the multiplicities of $N(E_i)$ and $\nu(E_i)-1$ of exceptional 
curves appearing in $f \circ \pi$ and $K_{X/{\mathbb C}^n}=
\pi^*(dx_1\wedge\ldots\wedge dx_n)$ which already appeared in the previous
section, we can now define and compute the topological Zeta-function.
As before, we first recall the definitions and properties (see e.g. 
\cite{DL}) and then continue with the computational aspects -- listing 
the corresponding {\sc Singular} commands in the continued example of an 
$A_4$ surface singularity.\\

\begin{definition}
Let $f \in {\mathbb C}[x_1,\ldots,x_n]$ be a non-zero polynomial defining a
hypersurface $V$ and let $\pi :X \longrightarrow {\mathbb C}^n$ be an embedded
resolution of $V$. Denote by $E_i$, $i \in I$, the irreducible components
of the divisor $\pi^{-1}(f^{-1}(0))$ . To fix notation, we define for each 
subset $J \subset I$
\begin{center}
$E_J:=\cap_{j \in J}E_j$   and   
$E_J^* := E_J \smallsetminus \cup_{j \notin J}E_{J \cup \{j\}}$ 
\end{center} 
and denote for each $j \in I$ the multiplicity of $E_j$ in the divisor 
of $f \circ \pi$ by $N(E_j)$. We further set $\nu(E_j)-1$ to be the 
multiplicity of $E_j$ in the divisor 
$K_{X/{\mathbb C}^n}=\pi^*(dx_1\wedge\ldots\wedge dx_n)$.
In this notation the topological Zeta-function of $f$ is 
      $$Z_{top}^{(d)}(f,s) := 
         \sum_{{{J \subset I {\rm s.th.}} \atop
                {d|N(E_j) \forall j \in J}}} \chi(E_J^*)
                \prod_{j \in J}(\nu(E_j) + N(E_j)s)^{-1} \in {\mathbb Q}(s).$$
Intersecting the $E_J^*$ with the preimage of zero in the above formula
leads to the local topological Zeta-function 
      $$Z_{top,0}^{(d)}(f,s) :=
         \sum_{{{J \subset I {\rm s.th.}} \atop
                {d|N(E_j) \forall j \in J}}} \chi(E_J^* \cap \pi^{-1}(0))
                \prod_{j \in J}(\nu(E_j) + N(E_j)s)^{-1} \in {\mathbb Q}(s).$$
(The local and global topological Zeta-function are independent of the choice
of embedded resolution of singularities of $V$.)
\end{definition}

Here, it is again important to observe that in the above context the
irreducible components are taken over ${\mathbb C}$, while practical 
calculations usually take place over ${\mathbb Q}$ and further passing to 
components taken over ${\mathbb C}$ is rather expensive. The following lemma 
shows that considering ${\mathbb Q}$-irreducible components already allows 
the computation of the $\zeta$-function:

\begin{lemma}
Let $D_l$ , $l \in L$, be the ${\mathbb Q}$-irreducible components of the
divisor $\pi^{-1}(f^{-1}(0))$. For each subset $J \subset L$ define $D_J$  
and $D_J^*$ as above. Then the topological (global and local) Zeta-function
can be computed by the above formulae using the $D_J$ and $D_J^*$ instead 
of the $E_J$ and $E_J^*$. 
\end{lemma}

As we already identified the exceptional divisors and computed the 
multiplicities $N(E_i)$ and $\nu(E_i)$, the only computational data 
missing is the Euler characteristic of the exceptional components in the
final charts. If the intersection of exceptional divisors is zero-dimensional,
this is just a matter of counting points using the identification of points
in different charts. For 1-dimensional intersections the Euler characteristic
can be computed using the geometric genus of the curve 
(using $\chi(C)=2-2g(C)$). Starting from dimension two on, this becomes more
subtle. 

\begin{example}
In our example which we have been treating throughout this article, we are
dealing with a surface in a three-dimensional ambient space. So the only
further Euler characteristics which need to be determined are those of the
exceptional divisors themselves. At the moment of birth of an exceptional
divisor, it will either be a ${\mathbb P}^2$ with Euler characteristic $3$
or a ${\mathbb P}^1 \times C$ (for a one-dimensional center $C$) leading to 
Euler characteristic $4 -4g(C)$. Under subsequent blow-ups the tracking of 
the changes to the Euler characteristic is then no difficult task. 

\begin{verbatim}
// compute the topological zeta-function for our isolated 
// surface singularity (global and local zeta-function coincide)
> zetaDL(L,1);           // global zeta-function
[1]:
   (s+6)/(5s2+11s+6)

> zetaDL(L,1,"local");   // local zeta-function
Local Case: Assuming that no (!) charts were dropped
during calculation of the resolution (option "A")
[1]:
   (s+6)/(5s2+11s+6)

// zetaDL also computes the characteristic polynomial
// of the monodromy, if additional parameter "A" is given
> zetaDL(L,1,"A");
Computing global zeta function
[1]:
   (s+6)/(5s2+11s+6)
[2]:
   (s4+s3+s2+s+1)
\end{verbatim}

\end{example}

\section{Desingularization in Experiments}

The previous section showed some examples in which desingularization was
a crucial step in the calculation of certain invariants and was hence used as 
a theoretical and practical tool. We now turn our interest to a different 
kind of settings: experiments on open questions which involve 
desingularization. Here we only sketch two such topics and the way one 
could experiment in the respective setting.\footnote{Neither of the two
topics should be seen as a suggestion for a short term research project!
Both questions, however, might gain new insights from someone playing 
around with such experiments just for a short while and stumbling
into examples which open up knew perspectives, insight or conjectures.}

\subsection{Bernstein-Sato polynomials}

In the early 1970s J.Bernstein \cite{Ber} and M.Sato \cite{Sato} independently 
defined an object in the theory of D-modules, which is nowadays called the
Bernstein-Sato polynomial. The roots of such polynomials have a close,
but still somewhat mysterious relation to the multiplicities of exceptional
divisors in a related desingularization. For briefly recalling the 
definition of Bernstein-Sato polynomials, we shall follow the article of 
Kashiwara \cite{Kash}, which also introduces this relation. After that 
we sketch what computer algebra tools are available in {\sc Singular} for 
experiments on this topic.

\begin{definition}
Let $f$ be an analytic function defined on some complex manifold $X$ of 
dimension $n$ and let ${\mathcal D}$ be the sheaf of differential operators 
of finite order on $X$. The polynomials in an additional variable $s$ 
satisfying
$$b(s) f^s \in {\mathcal D}[s] f^{s+1}$$
form an ideal. A generator of this ideal is called the Bernstein-Sato 
polynomial and denoted by $b_f(s)$.  
\end{definition}

Kashiwara then proves the rationality of the roots of the Bernstein-Sato
polynomial by using Hironaka's desingularization theorem in the following 
way:

\begin{theorem}[\cite{Kash}]
Let $f$ be as above and consider a blow-up $F:X' \longrightarrow X$ and
a new function $f'=f \circ F$. Then $b_f(s)$ is a divisor of 
$\prod_{k=0}^N b_{f'}(s+k)$ for a sufficiently large $N$.
\end{theorem}

More precisely, a principalization of the ideal generated by $f$ 
leads to a polynomial $f'=\prod_{i=1}{m} t_i^{r_i}$ with a local system 
of coordinates $t_1, \dots , t_m$. For $f'$ the Bernstein-Sato polynomial
is known to be
$$b_{f'}(s) = \prod_{i=1}^{m}\prod_{k=1}^{r_i} (r_i s + k).$$
Therefore the roots of the Bernstein-Sato polynomial are negative rational
numbers of the form $-\frac{k}{r_i}$ with exceptional mulitplicities $r_i$
and $1 \leq k \leq r_i$. Yet it is not clear at all whether there are
any rules or patterns which exceptional multiplicities
$r_i$ actually appear as denominators of roots of the Bernstein-Sato 
polynomial. 
(As a sideremark: If this were known, then such knowledge could be used 
to speed up computations of Bernstein-Sato polynomials whenever a 
desingularization is known.)\\

Today there are implemented algorithms for computing the Bernstein-Sato
polynomial of a given $f$ available: by Ucha and Castro-Jiminez \cite{UCJ},
Andres, Levandovskyy and Martin-Morales \cite{ALM} and by Berkesh and Leykin
\cite{BL}. All of these algorithms do not use desingularization techniques,
but rather rely on Gr\"obner-Bases and annihilator computations. The computed
objects, however, are also used (theoretically and in examples) in the context 
of multiplier ideals and thus have relation to invariants like the 
log-canonical threshold.\\

Given algorithmic ways to independently compute the Bernstein-Sato polynomial
and the exceptional multiplicities, one could now revisit the exploration of
the interplay between these and try to spot patterns to get a better 
understanding.  

\subsection{Positive Characteristic}

Desingularization in positive characteristic is one of the long standing,
central open problems in algebraic geometry. In dimensions up 
to three there is a positive answer (see e.g. \cite{CP1}, \cite{CP2}), but 
in the general case there are several different approaches (e.g. \cite{Ka}, 
\cite{Vi}, \cite{HW}) each of which has run into obstacles which are currently not
resolved.\\

About a decade ago, Hauser started studying the reasons why Hironaka's
approach of characteristic zero fails in positive characteristic \cite{Ha1}.
Among other findings, he singled out two central points which break down:
\begin{enumerate}
\item failure of maximal contact:\\
      Hypersurfaces of maximal contact are central to the descent in 
      ambient dimension which in turn is the key to finding the correct
      centers for blowing up. In positive characteristic, it is well 
      known that hypersurfaces of maximal contact need not exist;
      allowing hypersurfaces satisfying only slightly weaker conditions
      is one of the central steps in the approach of Hauser and Wagner
      for dimension 2 \cite{HW}. (For higher dimensions this definition
      requires a little bit more care \cite{AFK}.) 
\item increase of order of coefficient ideal:\\
      As the improvement of the singularities is measured by the decrease
      of the order and of orders for further auxilliary ideals constructed
      by means of descent in ambient dimension, it is crucial for the 
      proof of termination of resolution that these orders cannot increase
      under blowing up. Unfortunately this does no longer hold for the
      orders of the auxilliary ideals in positive characteristic as has 
      again been known since the 1970s. Hauser characterized the structure of 
      polynomials which can exhibit such behaviour in \cite{Ha1}.
\end{enumerate}

Although problems of desingularization in positive characteristic are known, 
this knowledge
seems to be not yet broad enough to provide sufficient feedback for
suitable modification of one of the approaches to overcome the respective
obstacles. Experiments could prove to be helpful to open up a new 
point of view. In particular, the approach of Hauser and Wagner for surfaces
is sufficiently close to the characteristic zero approach of Hironaka (and 
hence to algorithmic approaches like the one of Villamayor) to allow 
modification of an existing implementation to provide a tool for a 
structured search for examples with special properties also in higher 
dimensions. This has e.g. been pursued in \cite{AFK}.

\end{document}